\documentclass[11pt, english, draft]{article}
\usepackage{amssymb,amsmath}
\title{Several Quantitative Characterizations of Some Specific Groups}
\author{ {\bf A. Mohammadzadeh} \ and \ {\bf A. R. Moghaddamfar} \\[0.1cm]
{\em Faculty of Mathematics, K. N. Toosi
University of Technology,}\\[0.1cm]
{\em P. O. Box $16315$--$1618$, Tehran, Iran.}\\[0.1cm]
{\em E-mails}: {\tt moghadam@kntu.ac.ir} \ and \ {\tt
moghadam@ipm.ir}}

\newtheorem{theorem}{Theorem}[section]
\newtheorem{definition}[theorem]{Definition}
\newtheorem{corollary}[theorem]{Corollary}

\newtheorem{pro}[theorem]{Proposition}

\newtheorem{lm}[theorem]{Lemma}

\newtheorem{que}[theorem]{Question}
\begin{document}
\newcommand{\f}{\frac}
\newcommand{\sta}{\stackrel}
\maketitle
\begin{abstract}
\noindent Let $G$ be a finite group and let $\pi(G)=\{p_1, p_2,
\ldots, p_k\}$ be the set of prime divisors of $|G|$ for which
$p_1<p_2<\cdots<p_k$. The Gruenberg-Kegel graph of $G$, denoted
${\rm GK}(G)$, is defined as follows: its vertex set is $\pi(G)$
and two different vertices $p_i$ and $p_j$ are adjacent by an
edge if and only if $G$ contains an element of order $p_ip_j$.
The degree of a vertex $p_i$ in ${\rm GK}(G)$ is denoted by
$d_G(p_i)$ and the $k$-tuple $D(G)=\left(d_G(p_1), d_G(p_2),
\ldots, d_G(p_k)\right)$ is said to be the degree pattern of $G$.
Moreover, if $\omega \subseteq \pi(G)$ is the vertex set of a
connected component of ${\rm GK}(G)$, then the largest
$\omega$-number which divides $|G|$, is said to be an order
component of ${\rm GK}(G)$. We will say that the problem of
OD-characterization is solved for a finite group if we find the
number of pairwise non-isomorphic finite groups with the same
order and degree pattern as the group under study. The purpose of
this article is twofold. First, we completely solve the problem
of OD-characterization for every finite non-abelian simple group
with orders having prime divisors at most 29. In particular, we
show that there are exactly two non-isomorphic finite groups with
the same order and degree pattern as $U_4(2)$. Second, we prove
that there are exactly two non-isomorphic finite groups with the
same order components as $U_5(2)$.
\end{abstract}
\renewcommand{\baselinestretch}{1.1}
\def\thefootnote{ \ }
\footnotetext{{\em AMS subject Classification {\rm 2010}}:
20D05, 20D06, 20D08.\\[0.1cm]
}
\section{Introduction}
Throughout this article, all the groups under consideration are
{\em finite}, and simple groups are {\em non-abelian}. Given a
group $G$, the spectrum $\omega(G)$ of $G$ is the set of orders of
elements in $G$. Clearly, the spectrum $\omega(G)$ is closed and
partially ordered by the divisibility relation, and hence is
uniquely determined by the set $\mu(G)$ of its elements which are
maximal under the divisibility relation. If $n$ is a natural
number, then $\pi(n)$ denotes the set of all prime divisors of
$n$, in particular, we set $\pi(G)=\pi(|G|)$.

One of the most well-known graphs associated with $G$ is the
Gruenberg-Kegel graph (or prime graph) denoted by ${\rm GK}(G)$.
The vertex set of this graph is $\pi(G)$ and two distinct
vertices $p$ and $q$ are joined by an edge (abbreviated $p\sim
q$) if and only if $pq\in\omega(G)$. The number of connected
components of ${\rm GK}(G)$ is denoted by $s(G)$ and the sets of
vertices of connected components of ${\rm GK}(G)$ are denoted as
$\pi_i=\pi_i(G)( i =1, 2, \ldots, s(G))$. If $G$ is a group of
even order, then we put $2 \in \pi_1(G)$. The vertex sets of
connected components of all finite simple groups are obtained in
\cite{kondra} and \cite{wili}.

Given a group $G$, suppose that $\pi(G)=\{p_1, p_2, \ldots,
p_k\}$ in which $p_1 <p_2 < \cdots <p_k$. The {\em degree}
$d_G(p_i)$ of a vertex $p_i$ in the prime graph ${\rm GK}(G)$ is
the number of edges incident on $p_i$. We define
$D(G)=\left(d_G(p_1), d_G(p_2), \ldots, d_G(p_k )\right)$, and we
call this $k$-tuple the {\em degree pattern of} $G$. In addition,
we denote by $\mathcal{OD}(G)$ the set of pairwise non-isomorphic
finite groups with the same order and degree pattern as $G$, and
we put $h_{\rm OD}(G)=|\mathcal{OD}(G)|$. Clearly, there are only
finitely many isomorphism types of groups of order $|G|$, because
there are just finitely many ways that an $|G|\times |G|$
multiplication table can be filled in. Finally, for each group
$G$, it is clear that $1\leqslant h_{\rm OD}(G)<\infty$.

\begin{definition}
A group $G$ is called $k$-fold OD-characterizable if $h_{\rm
OD}(G)=k$. Usually, a $1$-fold OD-characterizable group is simply
called OD-characterizable, and it is called quasi
OD-characterizable if it is $k$-fold OD-characterizable for some
$k>1$.
\end{definition}

We will say that the OD-characterization problem is solved for a
group $G$, if the value of $h_{\rm OD}(G)$ is known. Studies in
recent years by several researchers have shown that many simple
groups are OD-characterizable. Some of these results are
summarized in the following table.
\begin{center}
\begin{tabular}{|l|l|c|l|}
\multicolumn{4}{c}{\textbf{Table 1.} Some OD-characterizable groups.} \\[0.1cm]
\hline $G$ & Condition on $G$ & $h_{\rm OD}(G)$ & References \\[0.1cm]
\hline
${\Bbb A}_n$ & $n=p, p+1, p+2$ ($p$ a prime, $p\geqslant 5$) & 1 & \cite{M.Z(2008)}, \cite{degree} \\[0.1cm]
& $5\leqslant n\leqslant 100, n\neq10$ & 1 & \cite{H.M(2010)}, \cite{kogani}, \cite{$S_6(3)$},\\[0.1cm]
&&& \cite{M.Z(2009)}, \cite{$A_{16}$} \\[0.1cm]
& $n=106, 112, 116, 134$ & 1 & \cite{chen(2010)}, \cite{chen(2013)}\\[0.1cm]
$L_2(q)$ & $q\neq 2, 3$& 1 & \cite{degree}, \cite{$L_2(q)$}\\[0.1cm]
$L_3(q)$ & $\left|\pi\left(\frac{q^2+q+1}{d}\right)\right|=1$, $d=(3, q-1)$ & 1& \cite{degree} \\[0.1cm]
$L_4(q)$ & $q\leqslant 17$ and $q=19, 23, 27, 29, 31, 32, 37$ & 1 & \cite{BAk2012}, \cite{four},  \cite{MAk.Rah} \\[0.1cm]
$L_n(2)$ & $n=p$ or $p+1$, $2^p-1$ is Mersenne prime & 1&\cite{MAk.Rah} \\[0.1cm]
$L_n(2)$ & $n=9, 10, 11$ & 1 & \cite{$L_9(2)$}, \cite{binary} \\[0.1cm]
$L_3(9)$ & & 1& \cite{$L_3(9)$}\\[0.1cm]
$L_6(3)$ & & 1& \cite{BAk(sub)}\\[0.1cm]
$U_3(q)$ & $\left|\pi\left(\frac{q^2-q+1}{d}\right)\right|=1$, $d=(3,q+1)$, $q>5$ & 1& \cite{degree} \\[0.1cm]
$U_4(q)$ & $q= 5 ,7$ & 1 &\cite{BAk(sub)}, \cite{MAk.Rah}\\[0.1cm]
$R(q)$ & $\left|\pi(q \pm \sqrt{3q}+1)\right|=1, q=3^{2m+1},
m\geqslant1$
& 1 & \cite{degree} \\[0.1cm]
${\rm Sz}(q)$ & $q=2^{2n+1}\geqslant 8$ & 1& \cite{degree} \\[0.1cm]
$O_5(q)\cong S_4(q)$ & $|\pi((q^2+1)/2)|=1$, $q\neq 3$ & 1& \cite{2-fold}\\[0.1cm]
$O_{2n+1}(q)\cong S_{2n}(q)$ & $n=2^m\geqslant 2, \ 2\mid q, \
|\pi(q^n+1)|=1$ &1&\cite{2-fold} \\[0.1cm]
&$(n, q)\neq(2, 2)$ &&\\[0.1cm]
$S_6(4)$ & & 1& \cite{$C_3(4)$} \\[0.1cm]
$G$ & $G$ is a sporadic group & 1& \cite{degree} \\[0.1cm]
$G$ & $|G|\leqslant 10^8 , G\neq {\Bbb A}_{10}, U_4(2)$ & 1& \cite{$10^8$}\\[0.1cm]
$G$ & $|\pi(G)|=4, G\neq {\Bbb A}_{10}$ & 1 & \cite{$K_4$}\\[0.1cm]
$G$ & $G$ is a simple with $\pi_1(G)=\{2\}$ & 1 & \cite{M.Z(2008)}\\[0.1cm]
\hline
\end{tabular}
\end{center}

In connection with the simple groups which are quasi
OD-characterizable, it was shown in \cite{2-fold},
\cite{connected} and \cite{degree} that:
\begin{center}
$\begin{array}{ll}
\mathcal{OD}(\mathbb{A}_{10})=\{\mathbb{A}_{10}, \ {\Bbb
Z}_3\times J_2\}, & \\[0.2cm]
 \mathcal{OD}(S_{6}(5))=\{ S_6(5), O_{7}(5)\}, &  \\[0.2cm]
\mathcal{OD}(S_{2m}(q))=\{S_{2m}(q), O_{2m+1}(q)\}, &
m=2^f\geqslant 2, \ \
\left|\pi\left(\frac{q^m+1}{2}\right)\right|=1,
\  \ q \ {\rm odd \ prime \ power},\\[0.2cm]
\mathcal{OD}(S_{2p}(3))=\{S_{2p}(3), O_{2p+1}(3)\}, &
\left|\pi\left(\frac{3^p-1}{2}\right)\right|=1, \ \ p \ {\rm odd \
prime}.  \\[0.2cm]
\end{array}$
\end{center}
In addition to the above results, it has been shown that in
\cite{moghadam} there exist many infinite families of alternating
and symmetric groups, $\{\mathbb{A}_n\}$ and $\{\mathbb{S}_n\}$,
which are quasi OD-characterizable, with $h_{\rm OD}(G)\geqslant
3$.

Here, we consider the simple groups $S$ such that
$\pi(S)\subseteq \pi(29!)$, and we denote the set of all these
simple groups by $\mathcal{S}_{\leqslant 29}$. Using the
classification of finite simple groups it is not hard to obtain a
full list of all groups in $\mathcal{S}_{\leqslant 29}$.
Actually, there are 110 such groups (see \cite[Table 4]{$S_6(3)$}
or \cite[Table 1]{za}). For convenience, the groups $S$ in
$\mathcal{S}_{\leqslant 29}$  and their orders are listed in
Table 2. The comparison between simple groups listed in Table 1
and the simple groups in $\mathcal{S}_{\leqslant 29}$ shows that
there are only 5 groups in $\mathcal{S}_{\leqslant 29}$, that is
$L_3(11)$, $U_4(2^3)$, ${^2E}_6(2)$, $S_4(17)$ and $U_4(17)$,
that the OD-characterization problem has not solved for them.
Therefore, one of the aims of this article is to prove these
groups are OD-characterizable.
\begin{theorem} \label{asli}
The simple groups $L_3(11)$, $U_4(2^3)$, ${^2E}_6(2)$, $S_4(17)$
and $U_4(17)$ are OD-characterizable.
\end{theorem}

We recall that a group $G$ is called a $2$-Frobenius group if
$G=ABC$, where $A$ and $AB$ are normal subgroups of $G$, $B$ is a
normal subgroup of $BC$, and $AB$ and $BC$ are Frobenius groups.
Zinov\`eva and V. D. Mazurov observed that the prime graph of a
2-Frobenius group is always disconnected, more precisely, it is
the union of two connected components each of which is a complete
graph \cite[Lemma 3(a) ]{zin-maz}. On the other hand, Mazurov
constructed a $2$-Frobenius group of the same order as the simple
group $U_4(2)$ (\cite{mtaiwaese, zin-maz}). In particular, this
shows that $h_{\rm OD }(U_4(2))\geqslant 2$ (see also
\cite{$10^8$}). In this article we also prove that the following
result.

\begin{theorem} \label{u4(2)} The simple group $U_4(2)$ is $2$-fold
OD-characterizable. In fact, there exists a unique $2$-Frobenius
group $F=(2^4\times 3^4):5:4$ with the same order and degree
pattern as $U_4(2)$, and so $\mathcal{OD}(U_4(2))=\{U_4(2), F\}$.
\end{theorem}

As an immediate consequence of Theorems \ref{asli},  \ref{u4(2)}
and the results in \cite{mtaiwaese, connected, degree}, we have
the following corollary.
\begin{corollary}
All simple groups in the class  $\mathcal{S}_{\leqslant 29}$,
other than ${\Bbb A}_{10}$, $S_6(3)$, $O_7(3)$ and $U_4(2)$, are
OD-characterizable. In addition, each of these groups is $2$-fold
OD-characterizable.
\end{corollary}

Given a group $G$, the order of $G$ can be expressed as a product
of some coprime natural numbers $m_i=m_i(G)$, $i=1, 2, \ldots,
s(G)$, with $\pi(m_i)=\pi_i$. The numbers $m_1, m_2, \ldots,
m_{s(G)}$ are called the {\it order components of} $G$. We set
$${\rm OC}(G)=\left\{m_1, m_2, \ldots,
m_{s(G)}\right\}.$$

In the similar manner, we denote by $\mathcal{OC}(G)$ the set of
isomorphism classes of finite groups with the same set ${\rm
OC}(G)$ of order components, and we put $h_{\rm
OC}(G)=|\mathcal{OC}(G)|$. Again, in terms of function $h_{\rm
OC}(\cdot)$, the groups $G$ are classified as follows:
\begin{definition} A group $G$ is called $k$-fold
OC-characterizable if $h_{\rm OC}(G)=k$. Usually, a $1$-fold
OC-characterizable group is simply called OC-characterizable, and
it is called quasi OC-characterizable if it is $k$-fold
OC-characterizable for some $k>1$.
\end{definition}

Obviously, if $p$ is a prime number, then $h_{\rm
OC}(\mathbb{Z}_p)=1$ while $h_{\rm OC}(\mathbb{Z}_{p^2})=h_{\rm
OC}(\mathbb{Z}_p\times \mathbb{Z}_p )=2$. Examples of
OC-characterizable groups are abundant (see for instance
\cite{iranmanesh-A-K}, \cite{iranmanesh}, \cite{A-Kho} and
\cite{Kkosravi-new}). Also, one family examples of simple groups
$S$ with $h_{\rm OC}(S)=2$ is given in \cite{A-Kho}, namely
$$\mathcal{OC}(O_{2n+1}(q))=\mathcal{OC}(S_{2n}(q))=\{O_{2n+1}(q), S_{2n}(q)\} \ \ \ \ \ \ (q \ {\rm odd}, \ n=2^m\geqslant 4).$$

\begin{center}
{\bf Table 2}. {\em Simple groups with orders having prime
divisors at most $29$ except alternating ones.}\\[0.3cm]

$\begin{array}{|ll|}\hline
S & |S|\\ \hline & \\[-0.26cm]
U_4(2) & 2^6\cdot 3^4\cdot 5 \\

L_2(7) &  2^3\cdot 3\cdot 7  \\

L_2(2^3) & 2^3\cdot 3^2\cdot 7  \\

U_3(3) & 2^5\cdot 3^3\cdot 7  \\

L_2(7^2) & 2^4\cdot 3\cdot 5^2\cdot 7^2 \\

U_3(5) & 2^4\cdot 3^2\cdot 5^3\cdot 7 \\

L_3(2^2) & 2^6\cdot 3^2\cdot 5\cdot 7   \\

J_2  & 2^7\cdot 3^3\cdot 5^2\cdot 7   \\

U_4(3) & 2^7\cdot 3^6\cdot 5\cdot 7 \\

S_4(7) & 2^8\cdot 3^2\cdot 5^2\cdot 7^4  \\

S_6(2) &  2^9\cdot 3^4\cdot 5\cdot 7  \\

O_8^+(2) & 2^{12}\cdot 3^5\cdot 5^2\cdot 7  \\

L_2(11) & 2^2\cdot 3\cdot 5\cdot 11  \\

M_{11} & 2^4\cdot 3^2\cdot 5\cdot 11  \\

M_{12} & 2^6\cdot 3^3\cdot 5\cdot 11  \\

U_5(2) & 2^{10}\cdot 3^5\cdot 5\cdot 11    \\

M_{12} & 2^6\cdot 3^3\cdot 5\cdot 11   \\

M^cL & 2^7\cdot 3^6\cdot 5^3\cdot 7\cdot 11  \\

HS & 2^{9}\cdot 3^2\cdot 5^3\cdot 7\cdot 11\\

U_6(2) & 2^{15}\cdot 3^6\cdot 5\cdot 7\cdot 11\\

L_3(3) &  2^{4}\cdot 3^3\cdot 13\\

L_2(5^2) & 2^3\cdot 3\cdot 5^2\cdot 13\\

U_3(2^2) & 2^{6}\cdot 3\cdot 5^2\cdot 13\\

S_4(5) & 2^6\cdot 3^2\cdot 5^4\cdot 13\\

L_4(3) &  2^{7}\cdot 3^6 \cdot 5\cdot 13\\

{{^2F}_4(2)}' & 2^{11}\cdot 3^3\cdot 5^2\cdot 13\\

L_2(13) &  2^2\cdot 3\cdot 7\cdot 13 \\

L_2(3^3) & 2^{2}\cdot 3^3\cdot 7\cdot 13  \\

G_2(3) & 2^{6}\cdot 3^6\cdot 7\cdot 13  \\

{^3D}_4(2) & 2^{12}\cdot 3^4\cdot 7^2\cdot 13  \\

{\rm Sz}(2^3) & 2^{6}\cdot 5\cdot 7\cdot 13 \\

L_2(2^6) & 2^{6}\cdot 3^2\cdot 5\cdot 7\cdot 13  \\

U_4(5) & 2^{7}\cdot 3^4\cdot 5^6\cdot 7\cdot 13 \\

L_3(3^2) & 2^{7}\cdot 3^6\cdot 5\cdot 7\cdot 13 \\

S_6(3) & 2^{9}\cdot 3^9\cdot 5\cdot 7\cdot 13 \\

O_7(3) & 2^{9}\cdot 3^9\cdot 5\cdot 7\cdot 13  \\

G_2(2^2) & 2^{12}\cdot 3^3\cdot 5^2\cdot 7\cdot 13  \\

S_4(2^3) & 2^{12}\cdot 3^4\cdot 5\cdot 7^2\cdot 13  \\

O_8^+(3) & 2^{12}\cdot 3^{12}\cdot 5^2\cdot 7\cdot 13  \\

L_5(3) & 2^9\cdot 3^{10}\cdot 5 \cdot {11}^2\cdot 13  \\

L_6(3) & 2^{11}\cdot 3^{15}\cdot 5\cdot 7\cdot {11}^2\cdot {13}^2\\

Suz & 2^{13}\cdot 3^7\cdot 5^2\cdot 7\cdot 11\cdot 13 \\

Fi_{22} & 2^{17}\cdot 3^{9}\cdot 5^2\cdot 7\cdot 11\cdot 13 \\
\hline
\end{array}$
$\begin{array}{|ll|}\hline
S & |S| \\
\hline & \\[-0.26cm]
L_2(17) & 2^4\cdot 3^2\cdot 17  \\

L_2(2^4) & 2^4\cdot 3\cdot 5\cdot 17  \\

S_4(2^2) &  2^8\cdot 3^2\cdot 5^2\cdot 17 \\

He & 2^{10}\cdot 3^3\cdot 5^2\cdot 7^3\cdot 17 \\

O_{8}^-(2) & 2^{12}\cdot 3^4\cdot 5\cdot 7\cdot 17  \\

L_4(2^2) & 2^{12}\cdot 3^4\cdot 5^2\cdot 7\cdot 17 \\

S_8(2) & 2^{16}\cdot 3^5\cdot 5^2\cdot 7\cdot 17 \\

U_4(2^2) & 2^{12}\cdot 3^2\cdot 5^3\cdot 13\cdot 17 \\

U_3(17) & 2^6\cdot 3^4\cdot 7\cdot 13\cdot 17^3 \\

O_{10}^-(2) & 2^{20}\cdot 3^6\cdot 5^2\cdot 7\cdot 11\cdot 17\\

L_2(13) &  2^2\cdot 3\cdot 7\cdot 13 \\

L_2(13^2) & 2^3\cdot 3\cdot 5\cdot 7\cdot 13^2\cdot 17 \\

S_{4}(13) & 2^6\cdot 3^2\cdot 5\cdot 7^2\cdot 13^4\cdot 17 \\

L_3(2^4) & 2^{12}\cdot 3^2\cdot 5^2\cdot 7\cdot 13\cdot 17\\

S_6(2^2) & 2^{18}\cdot 3^4\cdot 5^3\cdot 7\cdot 13\cdot 17\\

O_8^+(2^2) & 2^{24}\cdot 3^5\cdot 5^4\cdot 7\cdot 13\cdot 17^2\\

F_4(2) & 2^{24}\cdot 3^6\cdot 5^2\cdot 7^2\cdot 13\cdot 17 \\

{^2E}_6(2)&2^{36}\cdot 3^9\cdot 5^2\cdot 7^2\cdot 11\cdot13\cdot 17\cdot 19  \\

U_3(2^3)&2^9\cdot 3^4\cdot 7\cdot 19\\

U_4(2^3)&2^{18}\cdot 3^7\cdot 5\cdot 7^2\cdot 13\cdot 19\\

L_3(7)&2^5\cdot 3^2\cdot 7^3\cdot 19\\

L_4(7)&2^9\cdot 3^4\cdot 5^2\cdot 7^6\cdot 19\\

L_3(11)&2^4\cdot 3\cdot 5^2\cdot 7\cdot 11^3\cdot 19\\

L_2(19)&2^2\cdot 3^2\cdot 5\cdot 19\\

U_3(19)&2^5\cdot 3^2\cdot 5^2\cdot 7^3\cdot 19^3\\

J_1&2^3\cdot 3\cdot 5\cdot 7\cdot 11\cdot 19\\

J_3&2^7\cdot 3^5\cdot 5\cdot 17\cdot 19\\

F_5&2^{14}\cdot 3^6\cdot 5^6\cdot 7\cdot 11\cdot 19\\

L_2(23)&2^3\cdot 3\cdot 11\cdot 23\\

U_3(23)&2^7\cdot 3^2\cdot 11\cdot 13^2\cdot 23^3\\

M_{23}&2^7\cdot 3^2\cdot 5\cdot 7\cdot 11\cdot 23\\

M_{24}&2^{10}\cdot 3^3\cdot 5\cdot 7\cdot 11\cdot 23\\

Co_1&2^{21}\cdot 3^9\cdot 5^4\cdot 7^2\cdot 11\cdot 13\cdot23\\

Co_2&2^{18}\cdot 3^6\cdot 5^3\cdot 7\cdot 11\cdot 23\\

Co_3&2^{10}\cdot 3^7\cdot 5^3\cdot 7\cdot 11\cdot 23\\

Fi_{23}&2^{18}\cdot 3^{13}\cdot 5^2\cdot 7\cdot 11\cdot
13\cdot 17\cdot 23\\

U_4(17)&2^{11}\cdot 3^7\cdot 5\cdot 7\cdot 13\cdot 17^6\cdot29\\

S_4(17)&2^{10}\cdot 3^4\cdot 5\cdot 17^4\cdot 29\\

L_2(17^2)&2^5\cdot 3^2\cdot 5\cdot 17^2\cdot 29\\

L_2(29)&2^2\cdot 3\cdot 5\cdot 7\cdot 29\\

Ru & 2^{14}\cdot 3^3\cdot 5^3\cdot 7\cdot 13\cdot 29\\

Fi'_{24}& 2^{21}\cdot 3^{16}\cdot 5^2\cdot 7^3\cdot 11\cdot 13\cdot 17 \\
&  \cdot 23\cdot 29 \\
\hline
\end{array}$\end{center}

As the reader might have noticed, the values of the functions
$h_{\rm OD}$ and $h_{\rm OC}$ may be different. For example,
there are exactly two non-isomorphic groups of order
$1814400=2^7\cdot 3^4\cdot 5^2\cdot 7$ and degree pattern $(2, 3,
2, 1)$, they are $\mathbb{A}_{10}$ and $\mathbb{Z}_3\times J_2$,
and hence $h_{\rm OD}(\mathbb{A}_{10})=2$. However, since the
prime graph ${\rm GK}( \mathbb{A}_{10})$ is connected, ${\rm
OC}(\mathbb{A}_{10})=\{|\mathbb{A}_{10}|\}$, and so we obtain
$h_{\rm OC}(\mathbb{A}_{10})>\nu_{\rm a}(|\mathbb{A}_{10}|)=150$,
where $\nu_{\rm a}(m)$ denotes the number of types of abelian
groups of order $m$. Therefore, we have $h_{\rm
OD}(\mathbb{A}_{10})\neq h_{\rm OC}(\mathbb{A}_{10})$. The simple
group $U_5(2)$ is another example of this type. On the one hand,
we have $h_{\rm OD}(U_5(2))=1$ by Theorem 3.3 in \cite{$K_4$}. On
the other hand, there exists a $2$-Frobenius group $F$ such that
$|F|=|U_5(2)|$ (see \cite{mtaiwaese}) which implies that $h_{\rm
OC}(U_5(2))\geqslant 2$. Hence, $h_{\rm OD}(U_5(2))<h_{\rm
OC}(U_5(2))$. Finally, we show the following.
\begin{theorem} \label{u5(2)} The simple group $U_5(2)$ is $2$-fold
OC-characterizable. In fact, there exists a unique $2$-Frobenius
group $F=(2^{10}\times 3^5):11:5$ with the same order components
as $U_5(2)$, and so we have $\mathcal{OC}(U_5(2))=\{U_5(2), F\}$.
\end{theorem}

It is worth noting that the pair $\{U_5(2), (2^{10}\times
3^5):11:5\}$ is the first pair of a finite simple group and a
solvable group with the same order components. Note that these
groups have different prime graphs: the first connected component
of ${\rm GK}(U_5(2))$ is the path $2\sim 3\sim 5$, while the
first connected component of $(2^{10}\times 3^5):11:5$ is the
complete subgraph $2\sim 3\sim 5\sim 2$.

We introduce some more notation. Let $\Gamma$ be a simple graph.
An independent set of vertices in $\Gamma$ is a set of vertices
that are pairwise non-adjacent to each other in $\Gamma$. We
denote by $\alpha(\Gamma)$ the maximal number of vertices in
independent sets of $\Gamma$. Given a group $G$, we put
$t(G)=\alpha({\rm GK}(G))$. Moreover, for each prime $r\in
\pi(G)$, $t(r, G)$ denotes the maximal number of vertices in
independent sets of ${\rm GK}(G)$ containing $r$. Generally, our
notation for simple groups follows \cite{atlas}. Especially, the
alternating and symmetric group on $n$ letters are denoted by
${\Bbb A}_n$ and ${\Bbb S}_n$, respectively. We also denote by
${\rm Syl}_p(G)$ the set of all Sylow $p$-subgroups of $G$, where
$p\in\pi(G)$.

The sequel of this article is organized as follows: In Section 2,
we recall some basic results, especially, on the spectra of
certain finite simple groups. Section 3 is devoted to the proofs
of main results (Theorems \ref{asli}, \ref{u4(2)}, \ref{u5(2)}).
We conclude our article with  some open problems in Section 4.
\section{Preliminaries}
In this section, we consider some results which will be needed for
our further investigations.

\begin{lm}\label{vasi}{\rm (\cite{vasi})} Let $G$ be a finite group with
$t(G)\geqslant 3$ and $t(2, G)\geqslant 2$, and let $K$ be the
maximal normal solvable subgroup of $G$. Then the quotient group
$G/K$ is an almost simple group, i.e., there exists a non-abelian
simple group $P$ such that $P\leqslant G/K\leqslant {\rm Aut}(P)$.
\end{lm}

\begin{lm}{\rm \cite[Lemma 8]{luc}}\label{nonsolvable} Let $G$ be a finite group with $|\pi(G)|\geqslant 3$. If there exist
prime numbers $r, s, t \in\pi(G)$ such that $\{tr, ts,
rs\}\cap\omega(G)=\emptyset$, then $G$ is non-solvable.
\end{lm}

According to Table 4 in \cite{$S_6(3)$}, we have the following
result:
\begin{lm}
If $S\in\mathcal{S}_{\leqslant 29}$, then either ${\rm Out}(S)=1$
or $\pi({\rm Out}(S))\subseteq\{2, 3\}$.
\end{lm}

\begin{lm}\label{$L_2(q)$} {\rm (\cite{suzuki})}
Suppose that $q=p^n$, where $p$ is an odd prime. Then, we have
$$\mu(L_2(q))=\left\{p, \ \frac{q-1}{2}, \ \frac{q+1}{2}\right\}.$$
\end{lm}

\begin{lm}\label{$L_3(q)$} {\rm (\cite{m-comm})}
Suppose that $q=p^n$, where $p$ is an odd prime. Then, there
holds:
$$\mu(L_3(q))= \left\{
\begin{array}{ll}
\left\{q^2+q+1, \ q^2-1, \ p(q-1)\right\} & if \ \ \ q \not \equiv 1 \pmod{3},\\[0.3cm]
\left\{\frac{q^2+q+1}{3}, \ \frac{q^2-1}{3}, \ \frac{p(q-1)}{3}, \
q-1\right\} & if \ \ \ q \equiv 1 \pmod{3}.
\end{array}
\right.$$
\end{lm}

\begin{lm}\label{$L_4(q)$} {\rm (\cite{m-chen})}
Let $q$ be a power of prime $2$. Then, there holds:
$$\mu(U_4(q))=\left\{(q^2+1)(q-1), \ q^3+1, \ 2(q^2-1), \ 4(q+1)\right\}.$$
\end{lm}
\begin{lm}\label{$L_4(17)$} {\rm (\cite{zav-L4})}
Let $q$ be a power of an odd prime $p$. Denote $d={\rm gcd}(4,
q+1)$. Then $\mu(U_4(q))$ contains the following (and only the
following) numbers:
\begin{itemize}
\item[{\rm (i)}]  $\frac{q^4-1}{d(q+1)}$,
$\frac{q^3+1}{d}$, $\frac{p(q^2-1)}{d}$, $q^2-1$;
\item[{\rm (ii)}] $p(q+1)$, if and only if $d=4$;
\item[{\rm (iii)}] $9$, if and only if $p=3$.
\end{itemize}
\end{lm}
\begin{lm}\label{$S_4(q)$} {\rm (\cite{mazurov-2002})}
Let $q=p^n$, where $p>3$ is an odd prime. Then there holds:
$$\mu(S_4(q))=\left\{\frac{q^2+1}{2}, \ \frac{q^2-1}{2}, \
p(q+1), \ p(q-1)\right\}.$$
\end{lm}

Using Corollaries \ref{$L_3(q)$}, \ref{$L_4(q)$}, \ref{$L_4(17)$},
\ref{$S_4(q)$}, \cite[Table 4]{$S_6(3)$} and \cite{atlas} some
results are summarized in Table 3. In this table, we assume that
$s=|{\rm Out}(S)|$.

{\small
\begin{center}
\textbf{Table 3.} {\em Some simple groups in $\mathcal{S}_{\leqslant 29}$.}\\[0.5cm]
\begin{tabular}{|l|l|l|l|c|}
\hline $S$ & $|S|$& $\mu(S)$ & ${\rm D}(S)$ & $s$\\[0.1cm]
\hline $L_3(11)$ & $2^4\cdot3\cdot5^2\cdot7\cdot11^3\cdot19$& 110,
120, 133&(3, 2, 3, 1, 2, 1) &2\\[0.1cm]
$U_4(2^3)$ & $2^{18}\cdot 3^7\cdot 5\cdot 7^2\cdot 13\cdot 19$ &
36,
126, 455, 513 & (2, 3, 2, 4, 2, 1) & 6\\[0.1cm] $^2E_6(2)$ &
$2^{36}\cdot 3^9\cdot 5^2\cdot 7^2\cdot 11\cdot 13\cdot 17\cdot
19$ & 13, 16, \ldots, 22, 24, 28, 30, 33, 35 & (4, 4, 3, 3, 2,
0, 0, 0) & 6\\[0.1cm]
$S_4(17)$ & $2^{10}\cdot 3^4\cdot 5\cdot 17^4\cdot 29$ & 144, 145,
272, 306 & (2,
2, 1, 2, 1) & 2\\[0.1cm]
$U_4(17)$ & $2^{11}\cdot 3^7\cdot 5\cdot 7\cdot 13\cdot 17^6\cdot
29$ &
288, 2320, 2448, 2457 & (4, 4, 2, 2, 2, 2, 2) & 4\\[0.1cm] \hline
\end{tabular}
\end{center}}

The following proposition is taken from \cite{$10^8$}.

\begin{pro}{\rm (\cite{$10^8$})}\label{shi-zhang}
Let $M$ be a simple group whose order is less than $10^8$. If $G$
is a finite group with the same order and degree pattern as $M$,
then the following statements hold:
\begin{itemize}
\item[{\rm (a)}] If $M\neq A_{10}, U_4(2)$, then $G\cong M$;
\item[{\rm (b)}] If $M=A_{10}$, then $\mathcal{OD}(M)=\{A_{10}, J_2\times \mathbb{Z}_3\}$;
\item[{\rm (c)}] If $M=U_4(2)$, then $G$ is isomorphic to $M$ or a $2$-Frobenius group.
\end{itemize}
\end{pro}

In particular, Item (c) of Proposition \ref{shi-zhang} shows that
$h_{\rm OD}(U_4(2))\geqslant 2$. As we mentioned in the
Introduction, in fact, there is such a 2-Frobenius group (see
\cite{mtaiwaese}). Indeed, when we have a Frobenius group, say,
$F=K:C$ with abelian kernel $K$, and a faithful irreducible
$\mathbb{Z}_pF$-module $V$, then the semidirect product $VF$ is a
2-Frobenius group. Now, we consider the general linear groups
${\rm GL}(4,2)$ and ${\rm GL}(4,3)$. In ${\rm GL}(4,2)$ and also
in ${\rm GL}(4,3)$ there exists a Frobenius group $F=K:C$ of
order $20$ such that $K$ acts fixed-point-freely on corresponding
natural modules $V_1$  of dimension 4 over $\mathbb{F}_2$ and
$V_2$ of dimension 4 over $\mathbb{F}_3$. Now, we take
$(V_1\times V_2)\cdot E$ with the natural action of $E$ on direct
factors. Then we obtain a $2$-Frobenius group $(2^4\times
3^4):5:4$ with the same order as $U_4(2)$. Note that the prime
graphs of $U_4(2)$ and $(2^4\times 3^4):5:4$ coincide.

\section{Main Results}
In this section, we will prove Theorems \ref{asli}, \ref{u4(2)}
and \ref{u5(2)}. Before beginning the proof of Theorem
\ref{asli}, we draw the prime graphs of the groups $L_3(11)$,
$U_4(2^3)$, ${^2E}_6(2)$, $S_4(17)$ and $U_4(17)$ in Fig. 1 as
following:

{\small \setlength{\unitlength}{4mm}
\begin{picture}(0,0)(-39,-3)
\put(-31,-5){\circle*{0.35}}%
\put(-29,-5){\circle*{0.35}}%
\put(-27,-5){\circle*{0.35}}%
\put(-25,-5){\circle*{0.35}}%
\put(-30,-4){\circle*{0.35}}%
\put(-30,-6){\circle*{0.35}}%
\put(-28.6,-5.2){\footnotesize 2}%
\put(-31.8,-5.2){\footnotesize 5}%
\put(-27.25,-4.6){\footnotesize 7}%
\put(-25.3,-4.6){\footnotesize 19}%
\put(-30.1,-3.6){\footnotesize 3}%
\put(-30.4,-7){\footnotesize 11}%
\put(-31,-5){\line(1,0){2}}
\put(-31,-5){\line(1,1){1}}
\put(-31,-5){\line(1,-1){1}}
\put(-29,-5){\line(-1,1){1}}
\put(-29,-5){\line(-1,-1){1}}
\put(-27,-5){\line(1,0){2}}
\put(-30.5,-8.5){\footnotesize ${\rm GK}(L_3(11))$}

\put(-14,-5.5){\circle*{0.35}}%
\put(-12,-5.5){\circle*{0.35}}%
\put(-10,-5.5){\circle*{0.35}}%
\put(-8,-5.5){\circle*{0.35}}%
\put(-13,-4.5){\circle*{0.35}}%
\put(-11,-4.5){\circle*{0.35}}%
\put(-14.5,-6.5){\footnotesize 13}%
\put(-12.25,-6.5){\footnotesize 7}%
\put(-10.25,-6.5){\footnotesize 3}%
\put(-8.5,-6.5){\footnotesize 19}%
\put(-13.1,-4.1){\footnotesize 5}%
\put(-11.1,-4.1){\footnotesize 2}%
\put(-14,-5.5){\line(1,0){2}}
\put(-14,-5.5){\line(1,1){1}}
\put(-12,-5.5){\line(-1,1){1}}
\put(-12,-5.5){\line(1,0){2}}
\put(-12,-5.5){\line(1,1){1}}
\put(-10,-5.5){\line(1,0){2}}
\put(-10,-5.5){\line(-1,1){1}}
\put(-13.7,-8.5){\footnotesize ${\rm GK}(U_4(2^3))$}

\end{picture}}

\vspace{3.5cm}

{\small \setlength{\unitlength}{4mm}
\begin{picture}(0,0)(-42,-10)
\put(-39,-11){\circle*{0.35}}
\put(-35,-11){\circle*{0.35}}
\put(-37,-8.4){\circle*{0.35}}
\put(-37,-13.7){\circle*{0.35}}
\put(-37,-10){\circle*{0.35}}
\put(-33,-12){\circle*{0.35}}%
\put(-30,-12){\circle*{0.35}}%
\put(-31.5,-10){\circle*{0.35}}%
\put(-39.8,-11.2){\footnotesize 2}%
\put(-34.6,-11.2){\footnotesize 3}%
\put(-37.1,-8){\footnotesize 7}%
\put(-37.4,-14.7){\footnotesize 11}%
\put(-37.2,-10.9){\footnotesize 5}%
\put(-33.5,-13){\footnotesize 17}%
\put(-30.4,-13){\footnotesize 19}%
\put(-31.8,-9.5){\footnotesize 13}%
\put(-39,-11){\line(1,0){4}}
\put(-39,-11){\line(3,-4){2}}
\put(-39,-11){\line(3,4){2}}
\put(-39,-11){\line(2,1){2}}
\put(-35,-11){\line(-3,4){2}}
\put(-35,-11){\line(-2,1){2}}
\put(-35,-11){\line(-3,-4){2}}
\put(-37,-10){\line(0,1){1.65}}
\put(-37,-16){\footnotesize ${\rm GK}({^2E}_6(2))$}

\put(-24,-10){\circle*{0.35}}%
\put(-24,-12){\circle*{0.35}}%
\put(-23,-11){\circle*{0.35}}%
\put(-21,-11){\circle*{0.35}}%
\put(-19,-11){\circle*{0.35}}%
\put(-24.5,-9.7){\footnotesize 2}%
\put(-24.55,-12.7){\footnotesize 3}%
\put(-22.9,-10.7){\footnotesize 17}%
\put(-21.25,-10.5){\footnotesize 5}%
\put(-19.25,-10.5){\footnotesize 29}%
\put(-24,-10){\line(0,-1){2}}
\put(-24,-10){\line(1,-1){1}}
\put(-24,-12){\line(1,1){1}}
\put(-21,-11){\line(1,0){2}}
\put(-24,-16){\footnotesize ${\rm GK}(S_4(17))$}


\put(-14,-11){\circle*{0.35}}%
\put(-12,-11){\circle*{0.35}}%
\put(-10,-11){\circle*{0.35}}%
\put(-8,-11){\circle*{0.35}}%
\put(-13,-10){\circle*{0.35}}%
\put(-11,-10){\circle*{0.35}}%
\put(-9,-10){\circle*{0.35}}%
\put(-14.35,-12){\footnotesize 29}%
\put(-12.1,-12){\footnotesize 2}%
\put(-10.1,-12){\footnotesize 3}%
\put(-8.25,-12){\footnotesize 13}%
\put(-13.2,-9.65){\footnotesize 5}%
\put(-11.4,-9.65){\footnotesize 17}%
\put(-9.2,-9.65){\footnotesize 7}%
\put(-14,-11){\line(1,0){2}}
\put(-14,-11){\line(1,1){1}}
\put(-12,-11){\line(-1,1){1}}
\put(-12,-11){\line(1,0){2}}
\put(-12,-11){\line(1,1){1}}
\put(-10,-11){\line(1,0){2}}
\put(-10,-11){\line(-1,1){1}}
\put(-10,-11){\line(1,1){1}}
\put(-8,-11){\line(-1,1){1}}
\put(-13.5,-16){\footnotesize ${\rm GK}(U_4(17))$}

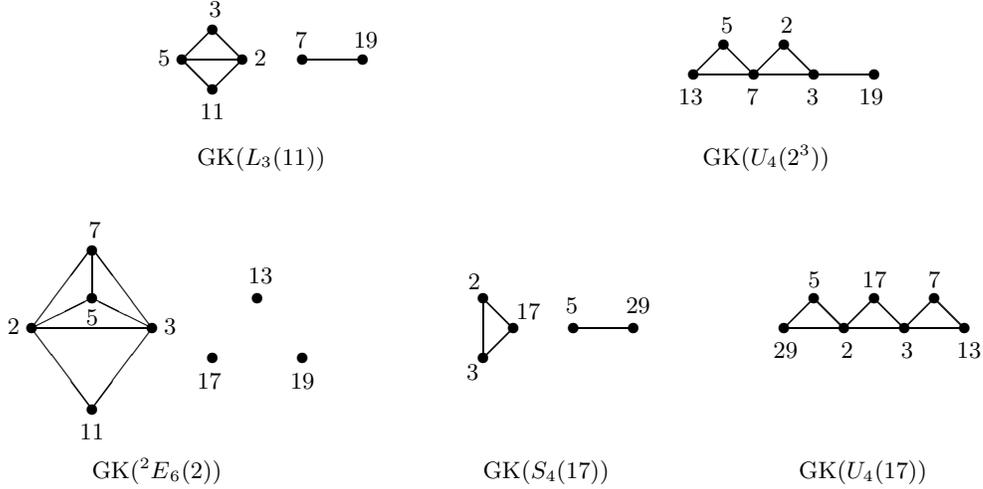
\put(-32,-18){\footnotesize {\bf Fig. 1.} The prime graph of some
simple groups.}
\end{picture}}
\vspace{4cm}

\noindent{\em Proof of Theorem $\ref{asli}$}. Let $S$ be one of
the following simple groups $L_3(11), U_4(2^3)$, ${^2E}_6(2)$ or
$U_4(17)$. Suppose that $G$ is a finite group such that $|G|=|S|$
and $D(G)=D(S)$. We have to prove that $G\cong S$. In all cases,
we will prove that $t(G)\geqslant 3$ and $t(2, G)\geqslant 2$.
Therefore, it follows from Lemma \ref{vasi} that there exists a
simple group $P$ such that $P\leqslant G/K\leqslant {\rm
Aut}(P)$, where $K$ is the maximal normal solvable subgroup of
$G$. In addition, we will prove that $P\cong S$, which implies
that $K=1$ and since $|G|=|S|$, $G$ is isomorphic to $S$, as
required. We handle every case singly.

\noindent (a) {\em $S=L_3(11)$}. Let $G$ be a finite group such
that
$$|G|=|S|=2^4\cdot 3\cdot 5^2\cdot 7\cdot 11^3\cdot 19 \ \
\mbox{and} \ \  D(G)=D(S)=(3, 2, 3, 1, 2, 1).$$

According to our hypothesis there are five possibilities for the
prime graph of $G$, as shown in Fig. 2. Here, $p_1, p_2\in\{2,
5\}$,  $p_3, p_4\in\{3, 11\}$, $p_5, p_6\in\{7, 19\}$.

\setlength{\unitlength}{4mm}
\begin{picture}(0,0)(-12,-4)
\put(-6.5,-7){\circle*{0.35}}%
\put(-4.5,-7){\circle*{0.35}}%
\put(-2,-7){\circle*{0.35}}%
\put(0,-7){\circle*{0.35}}%
\put(-5.5,-6){\circle*{0.35}}%
\put(-5.5,-8){\circle*{0.35}}%
\put(-7.25,-7.2){\footnotesize $2$}%
\put(-4.15,-7.2){\footnotesize $5$}%
\put(-2.2,-6.5){\footnotesize $7$}%
\put(-0.5,-6.5){\footnotesize $19$}%
\put(-5.6,-5.65){\footnotesize $3$}%
\put(-5.9,-8.9){\footnotesize $11$}%
\put(-6.5,-7){\line(1,0){2}}
\put(-6.5,-7){\line(1,1){1}}
\put(-6.5,-7){\line(1,-1){1}}
\put(-4.5,-7){\line(-1,1){1}}
\put(-4.5,-7){\line(-1,-1){1}}
\put(-2,-7){\line(1,0){2}}
\put(-4,-11){\footnotesize $(i)$}%

\put(3,-7){\circle*{0.35}}%
\put(5,-7){\circle*{0.35}}%
\put(7,-7){\circle*{0.35}}%
\put(9,-7){\circle*{0.35}}%
\put(11,-7){\circle*{0.35}}%
\put(8,-6){\circle*{0.35}}%
\put(2.75,-7.8){\footnotesize $p_5$}%
\put(4.75,-7.8){\footnotesize $p_4$}%
\put(6.75,-7.8){\footnotesize $p_2$}%
\put(8.75,-7.8){\footnotesize $p_1$}%
\put(10.75,-7.8){\footnotesize $p_6$}%
\put(7.75,-5.3){\footnotesize $p_3$}%
\put(3,-7){\line(1,0){2}}
\put(5,-7){\line(1,0){2}}
\put(7,-7){\line(1,1){1}}
\put(7,-7){\line(1,0){2}}
\put(9,-7){\line(-1,1){1}}
\put(9,-7){\line(1,0){2}}
\put(6,-11){\footnotesize $(ii)$}%

\put(14,-6){\circle*{0.35}}%
\put(14,-8){\circle*{0.35}}%
\put(16,-7){\circle*{0.35}}%
\put(18,-7){\circle*{0.35}}%
\put(20,-8){\circle*{0.35}}%
\put(20,-6){\circle*{0.35}}%
\put(13.7,-5.3){\footnotesize 7}%
\put(13.7,-9){\footnotesize 19}%
\put(15.8,-7.8){\footnotesize $p_1$}%
\put(17.5,-7.8){\footnotesize $p_2$}%
\put(19.6,-9){\footnotesize $11$}%
\put(19.85,-5.6){\footnotesize $3$}%
\put(16,-7){\line(-2,1){2}}
\put(16,-7){\line(-2,-1){2}}
\put(16,-7){\line(1,0){2}}
\put(18,-7){\line(2,1){2}}
\put(18,-7){\line(2,-1){2}}
\put(20,-8){\line(0,1){2}}
\put(16,-11){\footnotesize $(iii)$}%

\put(-1.5,-15){\circle*{0.35}}%
\put(0.5,-15){\circle*{0.35}}%
\put(2.5,-15){\circle*{0.35}}%
\put(4.5,-15){\circle*{0.35}}%
\put(0.5,-13.1){\circle*{0.35}}%
\put(2.5,-13.1){\circle*{0.35}}%
\put(-1.75,-15.9){\footnotesize 7}%
\put(0.25,-15.9){\footnotesize $p_1$}%
\put(2.25,-15.9){\footnotesize $p_2$}%
\put(4.25,-15.9){\footnotesize 19}%
\put(0,-12.4){\footnotesize $p_3$}%
\put(2.25,-12.4){\footnotesize $p_4$}%
\put(-1.5,-15){\line(1,0){6}}
\put(0.5,-15){\line(0,1){2}}
\put(2.5,-15){\line(0,1){2}}
\put(0.5,-13){\line(1,0){2}}
\put(0.5,-18){\footnotesize $(iv)$}%

\put(9.5,-15){\circle*{0.35}}%
\put(11.5,-15){\circle*{0.35}}%
\put(13.5,-15){\circle*{0.35}}%
\put(15.5,-15){\circle*{0.35}}%
\put(12.5,-14){\circle*{0.35}}%
\put(12.5,-16){\circle*{0.35}}%
\put(9.25,-14.6){\footnotesize 7}%
\put(10.7,-14.5){\footnotesize $p_1$}%
\put(13.5,-14.5){\footnotesize $p_2$}%
\put(15.25,-14.6){\footnotesize 19}%
\put(12.25,-13.6){\footnotesize $3$}%
\put(12.1,-17){\footnotesize $11$}%
\put(9.5,-15){\line(1,0){2}}
\put(11.5,-15){\line(1,1){1}}
\put(11.5,-15){\line(1,-1){1}}
\put(13.5,-15){\line(-1,1){1}}
\put(13.5,-15){\line(-1,-1){1}}
\put(13.5,-15){\line(1,0){2}}
\put(12,-18){\footnotesize $(v)$}%

\put(-2.5,-20){\footnotesize {\bf Fig. 2.} \ All possibilities
for the prime graph of $G$.}

\end{picture}
\vspace{7cm}

We now consider two subcases separately.
\begin{itemize}
\item[{\rm (a.1)}] {\em Assume first that ${\rm GK}(G)$ is disconnected}. In
this case, we immediately imply that ${\rm GK}(G)={\rm
GK}(L_3(11))$, and the hypothesis that $|G|=|L_3(11)|$ yields
${\rm OC}(G)={\rm OC}(L_3(11))$. Now, by the Main Theorem in
\cite{iranmanesh-A-K}, $G$ is isomorphic to $L_3(11)$, as
required.

\item[{\rm (a.2)}] {\em Assume next that ${\rm GK}(G)$ is connected}. In this
case $7\nsim 19$ in ${\rm GK}(G)$. Since $\{7, 19, p_3\}$ is an
independent set, $t(G)\geqslant 3$, and so by Lemma
\ref{nonsolvable}, $G$ is a non-solvable group. Moreover, since
$d_G(2)=3$ and $|\pi(G)|=6$, $t(2,G)\geqslant 2$. Thus by Lemma
\ref{vasi} there exists a simple group $P$ such that $P\leqslant
G/K\leqslant {\rm Aut}(P)$, where $K$ is the maximal normal
solvable subgroup of $G$. We claim that $K$ is a $\{7, 11,
19\}'$-group. We first show that $K$ is a $\{7, 19\}'$-group. If
$\{7, 19\}\subseteq \pi(K)$, then a Hall $\{7, 19\}$-subgroup of
$K$ is an abelian group. Hence $7\sim 19$ in ${\rm GK}(K)$, and
so in ${\rm GK}(G)$, that is a contradiction. Let $\{r, s\}=\{7,
19\}$. Now assume that $r\in\pi(K)$ and $s\notin\pi(K)$. Let $T\in
{\rm Syl}_{r}(K)$. By Frattini argument $G=KN_G(T)$. Therefore,
the normalizer $N_G(T)$ contains an element of order $s$, say $x$.
Now, $T\langle x\rangle$ is an abelian subgroup of $G$, so it
leads to a contradiction as before.

Finally, suppose that $11\in\pi(K)$ and $T\in {\rm Syl}_{11}(K)$.
Then, $G=KN_G(T)$, by Frattini argument. Evidently, $N_G(T)$
contains some elements of order 7 and 19, that we respectively
denote by $u$ and $v$. Now, $T\langle u\rangle$ and $T\langle
v\rangle$ are nilpotent subgroups of $G$, of orders $11^3\cdot 7$
and $11^3\cdot 19$, respectively, which implies that $7 \sim 11
\sim 19$, a contradiction. Since $K$ and ${\rm Out}(P)$ are $\{7,
11, 19 \}'$-groups, $|P|$ is divisible by $7\cdot11^3\cdot19$.
Considering the orders of simple groups in
$\mathcal{S}_{\leqslant 29}$, we conclude that $P$ is isomorphic
to $L_3(11)$, and so $K=1$ and since $|G|=|L_3(11)|$, $G$ is
isomorphic to $L_3(11)$. But then ${\rm GK}(G)={\rm GK}(L_3(11))$
is disconnected, which is impossible.
\end{itemize}

\noindent (b) {\em $S=U_4(2^3)$}. Assume that $G$ is a finite
group such that
$$|G|=|S|=2^{18}\cdot 3^7\cdot 5\cdot 7^2\cdot 13\cdot 19 \ \ \ {\rm
and} \ \ \ D(G)=D(S)=(2, 3, 2, 4, 2, 1).$$ So, the prime graph of
$G$ is one of the following graphs as shown in Fig. 3. Here $p_1,
p_2, p_3\in\{2, 5, 13\}$.

\setlength{\unitlength}{4mm}
\begin{picture}(0,0)(-24,-5.5)
\put(-22,-8){\circle*{0.35}}%
\put(-20,-8){\circle*{0.35}}%
\put(-18,-8){\circle*{0.35}}%
\put(-16,-8){\circle*{0.35}}%
\put(-21,-7){\circle*{0.35}}%
\put(-19,-7){\circle*{0.35}}%
\put(-22.25,-9){\footnotesize $p_1$}%
\put(-20.2,-9.1){\footnotesize 7}%
\put(-18.2,-9.1){\footnotesize 3}%
\put(-16.25,-9.1){\footnotesize 19}%
\put(-21.25,-6.5){\footnotesize $p_2$}%
\put(-19.25,-6.5){\footnotesize $p_3$}%
\put(-22,-8){\line(1,0){2}}
\put(-22,-8){\line(1,1){1}}
\put(-20,-8){\line(-1,1){1}}
\put(-20,-8){\line(1,1){1}}
\put(-20,-8){\line(1,0){2}}
\put(-18,-8){\line(1,0){2}}
\put(-18,-8){\line(-1,1){1}}
\put(-20,-11.8){\footnotesize $(i)$}%


\put(-12,-8){\circle*{0.35}}%
\put(-10,-8){\circle*{0.35}}%
\put(-8,-8){\circle*{0.35}}%
\put(-6,-8){\circle*{0.35}}%
\put(-11,-7){\circle*{0.35}}%
\put(-11,-9){\circle*{0.35}}%
\put(-12.3,-8.9){\footnotesize 3}%
\put(-9.9,-8.9){\footnotesize 7}%
\put(-8.25,-8.7){\footnotesize $p_3$}%
\put(-6.25,-8.9){\footnotesize 19}%
\put(-11.25,-6.5){\footnotesize $p_1$}%
\put(-11.25,-9.8){\footnotesize $p_2$}%
\put(-12,-8){\line(1,0){2}}
\put(-12,-8){\line(1,1){1}}
\put(-12,-8){\line(1,-1){1}}
\put(-10,-8){\line(1,0){2}}
\put(-10,-8){\line(-1,1){1}}
\put(-10,-8){\line(-1,-1){1}}
\put(-8,-8){\line(1,0){2}}
\put(-9.5,-11.8){\footnotesize $(ii)$}%


\put(-1.5,-8){\circle*{0.35}}%
\put(-1.5,-9.5){\circle*{0.35}}%
\put(0.5,-8){\circle*{0.35}}
\put(0.5,-9.5){\circle*{0.35}}%
\put(2.5,-8){\circle*{0.35}}%
\put(-0.5,-7){\circle*{0.35}}%
\put(-1.8,-7.7){\footnotesize 3}%
\put(-2.5,-10){\footnotesize $p_2$}%
\put(0.5,-7.7){\footnotesize 7}%
\put(1,-10){\footnotesize $p_3$}%
\put(2.25,-7.7){\footnotesize 19}%
\put(-0.75,-6.3){\footnotesize $p_1$}%
\put(-1.5,-8){\line(1,0){2}}
\put(-1.5,-8){\line(0,-1){1.5}}
\put(-1.5,-8){\line(1,1){1}}
\put(0.5,-8){\line(1,0){2}}
\put(0.5,-8){\line(-1,1){1}}
\put(0.5,-8){\line(0,-1){1.5}}
\put(-1.5,-9.5){\line(1,0){2}}
\put(-1,-11.8){\footnotesize $(iii)$}%


\put(7,-8){\circle*{0.35}}%
\put(9,-8){\circle*{0.35}}%
\put(11,-8){\circle*{0.35}}%
\put(13,-8){\circle*{0.35}}%
\put(9,-6.5){\circle*{0.35}}%
\put(9,-9.5){\circle*{0.35}}%
\put(6.75,-8.9){\footnotesize 3}%
\put(8.8,-8.9){\footnotesize $5$}%
\put(11.1,-8.9){\footnotesize 7}%
\put(12.75,-8.9){\footnotesize 19}%
\put(8.75,-6){\footnotesize $2$}%
\put(8.6,-10.6){\footnotesize $13$}%
\put(7,-8){\line(1,0){2}}
\put(7,-8){\line(4,3){2}}
\put(7,-8){\line(4,-3){2}}
\put(9,-8){\line(1,0){2}}
\put(11,-8){\line(1,0){2}}
\put(11,-8){\line(-4,3){2}}
\put(11,-8){\line(-4,-3){2}}
\put(9.2,-11.8){\footnotesize $(iv)$}%

\put(-14,-14){\footnotesize {\bf Fig. 3.} \ All possibilities for
the prime graph of $G$.}
\end{picture}
\vspace{4cm}

In what follows, we will consider two subcases separately.

\begin{itemize}
\item[{\rm (b.1)}] {\em First, suppose that ${\rm GK}(G)$ is one of the graphs
$(i), (iii)$ or $(iv)$}. Note that in each case $13\nsim 19$ in
${\rm GK}(G)$ and $t(G)\geqslant 3$. Now, it follows from Lemma
\ref{nonsolvable} that $G$ is a non-solvable group. Moreover,
since $d_G(2)=2$ and $|\pi(G)|=6$, $t(2, G)\geqslant 2$. Thus by
Lemma \ref{vasi}, there exists a simple group $P$ such that
$P\leqslant G/K\leqslant {\rm Aut}(P)$, where $K$ is the maximal
normal solvable subgroup of $G$. As the previous case, one can
show that $K$ is a $\{13, 19\}'$-group. Since $K$ and ${\rm
Out}(P)$ are $\{13, 19\}'$-groups, thus $|P|$ is divisible by
$13\cdot 19$. Considering the orders of simple groups in
$\mathcal{S}_{\leqslant 29}$ yields $P$ is isomorphic to
$U_4(8)$, and so $K=1$ and $G$ is isomorphic to $U_4(8)$, because
$|G|=|U_4(8)|$. Therefore, the prime graph of $G$ and the graph
$(i)$ coincide, and in other cases we get a contradiction.

\item[{\rm (b.2)}] {\em Next, suppose that ${\rm GK}(G)$ is the graph $(ii)$}.
In this case, $7$ is not adjacent to $19$ in ${\rm GK}(G)$. Since
$\{p_1, p_2, 19\}$ is an independent set, $t(G)\geqslant 3$ and by
Lemma \ref{nonsolvable}, $G$ is a non-solvable group. Moreover,
since $d_G(2)=2$ and $|\pi(G)|=6$, $t(2, G)\geqslant 2$. Thus by
Lemma \ref{vasi} there exists a simple group $P$ such that
$P\leqslant G/K\leqslant {\rm Aut}(P)$, where $K$ is the maximal
normal solvable subgroup of $G$. Using similar arguments to those
in the previous case, one can show that $K$ is a $\{7,
19\}'$-group and $G$ is isomorphic to $U_4(8)$. But then $3$ is
adjacent to $19$ in ${\rm GK}(G)$, which is a contradiction.
\end{itemize}

(c) {\em $S={^2E}_6(2)$}. Assume that $G$ is a finite group such
that
$$|G|=|S|=2^{36}\cdot 3^9\cdot 5^2\cdot 7^2\cdot 11\cdot 13\cdot 17\cdot
19 \ \ \ \mbox{and} \ \ \ D(G)=D(S)=(4, 4, 3, 3, 2, 0, 0, 0).$$
Then, the prime graphs of $G$ and ${^2E}_6(2)$ coincide, and the
hypothesis that $|G|=|{^2E}_6(2)|$ yields ${\rm OC}(G)={\rm
OC}({^2E}_6(2))$. Now, by \cite{Kkosravi-new}, $G$ is isomorphic
to ${^2E}_6(2)$, as required.

(d) {\em $S=U_4(17)$}. Assume that $G$ is a finite group such that
$$|G|=|S|=2^{11}\cdot 3^7\cdot 5\cdot 7\cdot 13\cdot 17^6\cdot 29 \ \ \
\mbox{and} \ \ \ D(G)=D(S)=(4, 4, 2, 2, 2, 2, 2).$$ According to
our hypothesis there are four possibilities for the prime graph of
$G$, as shown in Fig. 4. Here $p_1, p_2, p_3, p_4, p_5\in\{5, 7,
13, 17, 29\}$.

\setlength{\unitlength}{4mm}
\begin{picture}(0,0)(-32,-4)
\put(-31,-8.5){\circle*{0.35}}%
\put(-29,-8.5){\circle*{0.35}}%
\put(-27,-8.5){\circle*{0.35}}%
\put(-25,-8.5){\circle*{0.35}}%
\put(-30,-7.5){\circle*{0.35}}%
\put(-28,-7.5){\circle*{0.35}}%
\put(-26,-7.5){\circle*{0.35}}%
\put(-31.25,-9.2){\footnotesize $p_4$}%
\put(-29.2,-9.4){\footnotesize 2}%
\put(-27.2,-9.4){\footnotesize 3}%
\put(-25.25,-9.4){\footnotesize $p_5$}%
\put(-30.25,-7){\footnotesize $p_1$}%
\put(-28.25,-7){\footnotesize $p_2$}%
\put(-26.25,-7){\footnotesize $p_3$}%
\put(-31,-8.5){\line(1,0){2}}
\put(-31,-8.5){\line(1,1){1}}
\put(-29,-8.5){\line(-1,1){1}}
\put(-29,-8.5){\line(1,1){1}}
\put(-29,-8.5){\line(1,0){2}}
\put(-27,-8.5){\line(1,0){2}}
\put(-27,-8.5){\line(-1,1){1}}
\put(-27,-8.5){\line(1,1){1}}
\put(-25,-8.5){\line(-1,1){1}}
\put(-28.5,-13){\footnotesize $(i)$}%

\put(-20,-8){\circle*{0.35}}%
\put(-18,-8){\circle*{0.35}}%
\put(-16,-8){\circle*{0.35}}%
\put(-18,-6.5){\circle*{0.35}}%
\put(-18,-9.5){\circle*{0.35}}%
\put(-20,-10.5){\circle*{0.35}}%
\put(-16,-10.5){\circle*{0.35}}%
\put(-20.3,-7.7){\footnotesize 2}%
\put(-18.25,-7.5){\footnotesize $p_2$}%
\put(-16.1,-7.7){\footnotesize 3}%
\put(-18.25,-6){\footnotesize $p_1$}%
\put(-18.25,-8.8){\footnotesize $p_3$}%
\put(-20.25,-11.3){\footnotesize $p_4$}%
\put(-16.25,-11.3){\footnotesize $p_5$}%
\put(-20,-8){\line(1,0){4}}
\put(-20,-8){\line(4,3){2}}
\put(-20,-8){\line(4,-3){2}}
\put(-20,-8){\line(0,-1){2.5}}
\put(-16,-8){\line(-4,3){2}}
\put(-16,-8){\line(-4,-3){2}}
\put(-16,-8){\line(0,-1){2.5}}
\put(-20,-10.5){\line(1,0){4}}
\put(-18.5,-13){\footnotesize $(ii)$}%

\put(-11,-9){\circle*{0.35}}%
\put(-9,-9){\circle*{0.35}}%
\put(-7,-9){\circle*{0.35}}%
\put(-5,-9){\circle*{0.35}}%
\put(-9,-6.4){\circle*{0.35}}%
\put(-7,-6.4){\circle*{0.35}}%
\put(-8,-10){\circle*{0.35}}%
\put(-11.25,-9.9){\footnotesize $p_2$}%
\put(-9.5,-9.9){\footnotesize 2}%
\put(-7,-9.9){\footnotesize 3}%
\put(-5.25,-9.9){\footnotesize $p_3$}%
\put(-9.5,-6){\footnotesize $p_5$}%
\put(-7.25,-6){\footnotesize $p_4$}%
\put(-8.25,-10.8){\footnotesize $p_1$}%
\put(-11,-9){\line(1,0){2}}
\put(-11,-9){\line(3,4){2}}
\put(-9,-9){\line(1,0){2}}
\put(-9,-9){\line(3,4){2}}
\put(-9,-9){\line(1,-1){1}}
\put(-7,-9){\line(1,0){2}}
\put(-7,-9){\line(-3,4){2}}
\put(-7,-9){\line(-1,-1){1}}
\put(-5,-9){\line(-3,4){2}}
\put(-8.6,-13){\footnotesize $(iii)$}%

\put(-0.5,-9){\circle*{0.35}}%
\put(1.5,-9){\circle*{0.35}}%
\put(3.5,-9){\circle*{0.35}}%
\put(5.5,-9){\circle*{0.35}}%
\put(2.5,-4.5){\circle*{0.35}}%
\put(2.5,-8){\circle*{0.35}}%
\put(2.5,-10){\circle*{0.35}}%
\put(-0.75,-9.9){\footnotesize $p_4$}%
\put(1.2,-9.9){\footnotesize 2}%
\put(3.5,-9.9){\footnotesize 3}%
\put(5.25,-9.9){\footnotesize $p_5$}%
\put(2.25,-4){\footnotesize $p_1$}%
\put(2.25,-7.4){\footnotesize $p_2$}%
\put(2.25,-10.8){\footnotesize $p_3$}%
\put(-0.5,-9){\line(2,3){3}}
\put(1.5,-9){\line(-1,0){2}}
\put(1.5,-9){\line(1,0){2}}
\put(1.5,-9){\line(1,1){1}}
\put(1.5,-9){\line(1,-1){1}}
\put(3.5,-9){\line(1,0){2}}
\put(3.5,-9){\line(-1,1){1}}
\put(3.5,-9){\line(-1,-1){1}}
\put(5.5,-9){\line(-2,3){3}}
\put(2,-13){\footnotesize $(iv)$}%

\put(-22,-15){\footnotesize {\bf Fig. 4.} \ All possibilities for
the prime graph of $G$.}
\end{picture}
\vspace{5cm}

\noindent In all cases, $\{p_1, p_2, p_3\}$ is an independent set,
and hence $t(G)\geqslant 3$. Moreover, since $d_G(2)=4$ and
$|\pi(G)|=7$, $t(2, G)\geqslant 2$. Now, from Lemma \ref{vasi}
there exists a simple group $P$ such that $P\leqslant
G/K\leqslant {\rm Aut}(P)$. We claim now that $K$ is a $\{2,
3\}$-group. In fact, if there exists $p_i\in\pi(K)$, for some $i$,
then with similar arguments as before, we can verify that for each
$j\neq i$, $p_i\sim p_j$ in ${\rm GK}(G)$, and this contradicts
the fact that $d_G(p_i)=2$. Hence $K$ and ${\rm Out}(P)$ are $\{2,
3\}$-groups, thus $|P|$ is divisible by $5\cdot 7\cdot 13\cdot
17^6\cdot 29$. Again, considering the orders of simple groups in
$\mathcal{S}_{\leqslant 29}$ yields $P$ is isomorphic to
$U_4(17)$, and so $K=1$ and $G$ is isomorphic to $U_4(17)$,
because $|G|=|U_4(17)|$.

(e) {\em $S=S_4(17)$}. Assume that $G$ is a finite group such that
$$|G|=|S|=2^{10}\cdot 3^4\cdot 5\cdot 17^4\cdot 29 \ \ \ \mbox{and} \ \ \ D(G)=D(S)=(2, 2, 1, 2,
1).$$ Under these conditions, there are two possibilities for the
prime graph of $G$, as shown in Fig. 5. Here $p_1, p_2,
p_3\in\{2, 3, 17\}$.

\setlength{\unitlength}{4mm}
\begin{picture}(0,0)(-40,-6)
\put(-28.5,-7){\circle*{0.35}}%
\put(-28.5,-9){\circle*{0.35}}%
\put(-26.5,-8){\circle*{0.35}}%
\put(-23.5,-8){\circle*{0.35}}%
\put(-21.5,-8){\circle*{0.35}}%
\put(-28.7,-6.6){\footnotesize 2}%
\put(-28.75,-10){\footnotesize 3}%
\put(-26,-8.2){\footnotesize 17}%
\put(-23.7,-7.5){\footnotesize 5}%
\put(-21.75,-7.5){\footnotesize 29}%
\put(-28.5,-7){\line(0,-1){2}}
\put(-28.5,-7){\line(2,-1){2}}
\put(-28.5,-9){\line(2,1){2}}
\put(-23.5,-8){\line(1,0){2}}
\put(-26,-11){\footnotesize $(i)$}%

\put(-18,-8.5){\circle*{0.35}}%
\put(-17,-7.5){\circle*{0.35}}%
\put(-16,-8.5){\circle*{0.35}}%
\put(-15,-7.5){\circle*{0.35}}%
\put(-14,-8.5){\circle*{0.35}}%
\put(-18.2,-9.4){\footnotesize 5}%
\put(-17.5,-7){\footnotesize $p_1$}%
\put(-16.25,-9.4){\footnotesize $p_2$}%
\put(-15,-7){\footnotesize $p_3$}%
\put(-14.5,-9.4){\footnotesize 29}%
\put(-17,-7.5){\line(-1,-1){1}}
\put(-17,-7.5){\line(1,-1){1}}
\put(-15,-7.5){\line(-1,-1){1}}
\put(-15,-7.5){\line(1,-1){1}}
\put(-16.5,-11){\footnotesize $(ii)$}%


\put(-30,-13){\footnotesize {\bf Fig. 5.} \ All possibilities for
the prime graph of $G$.}
\end{picture}
\vspace{3.5cm}

\noindent We now consider two cases separately, depending on
${\rm GK}(G)$ is connected or disconnected.
\begin{itemize}
\item[{\rm (2.1)}] {\em Assume first that ${\rm GK}(G)$ is connected}.
Since $\{5, p_2, 29\}$ is an independent set, $t(G)\geqslant 3$.
Moreover, since $d_G(2)=2$ and $|\pi(G)|=5$, $t(2, G)\geqslant 2$.
Thus by Lemma \ref{vasi} there exists a simple group $P$ such
that $P\leqslant G/K\leqslant {\rm Aut}(P)$. We shall treat the
cases $17\nsim 29$ and $17\sim 29$ in ${\rm GK}(G)$, separately.

$(2.1.a)$ First we consider the case where $17\nsim 29$ in ${\rm
GK}(G)$. In this case as before, one can show that $K$ is a
$\{17, 29\}'$-group. Since $K$ and ${\rm Out}(P)$ are $\{17,
29\}'$-groups, thus $|P|$ is divisible by $17^4\cdot 29$.
Considering the orders of simple groups in
$\mathcal{S}_{\leqslant 29}$ yields $P$ is isomorphic to
$S_4(17)$, and so $K=1$ and $G$ is isomorphic to $S_4(17)$,
because $|G|=|S_4(17)|$. But then ${\rm GK}(G)= {\rm GK}(S_4(17))$
is disconnected, which is impossible.

$(2.1.b)$ Next we discuss the case where $17\sim 29$ in ${\rm
GK}(G)$. An argument similar to that in the above paragraphs
shows that $K$ is a $\{3, 29\}'$-group. Since $K$ and ${\rm
Out}(P)$ are $\{3, 29\}'$-groups, thus $|P|$ is divisible by
$3^4\cdot 29$. Considering the orders of simple groups in
$\mathcal{S}_{\leqslant 29}$ yields $P$ is isomorphic to
$S_4(17)$, and so $K=1$ and $G$ is isomorphic to $S_4(17)$,
because $|G|=|S_4(17)|$. But then, ${\rm GK}(G)= {\rm
GK}(S_4(17))$ is disconnected, which is impossible.

\item[{\rm (2.2)}] {\em Assume next that ${\rm GK}(G)$ is disconnected}. In
this case, it is easy to see that the prime graphs of $G$ and
$S_4(17)$ coincide. Now, by the main theorem in
\cite{iranmanesh}, $G$ is isomorphic to $S_4(17)$.
\end{itemize}
This completes the proof of Theorem $\ref{asli}$. $\Box$

\vspace{0.3cm}

\noindent{\em Proof of Theorem $\ref{u4(2)}$}. Let $G$ be a finite
group satisfying
$$(1) \ |G|=|U_4(2)|=2^6\cdot 3^4\cdot 5, \ \ \ {\rm and} \ \ \  (2) \ D(G)=D(U_4(2))=(1, 1, 0).$$
By Proposition \ref{shi-zhang}, $G$ is isomorphic to $U_4(2)$ or a
$2$-Frobenius group. First of all, it should be noted that the
existence of a $2$-Frobenius group satisfying conditions $(1)$ and
$(2)$ is guaranteed by \cite{mtaiwaese}. To prove uniqueness, we
note that any such group will be a subdirect product of
$2$-Frobenius groups of orders $2^4\cdot 5\cdot 4$ and $3^4\cdot
5\cdot 4$. As a matter of fact, since $4$ is the order of $2$
modulo $5$, $4$ is the smallest dimension of an irreducible
module for $\mathbb{Z}_5$ over $\mathbb{F}_2$, so there is a
unique Frobenius group of order $2^4\cdot 5$ and its kernel is
elementary abelian. Actually, this is a subgroup of the
$1$-dimensional affine group over $\mathbb{F}_{2^4}$ which is
denoted by ${\rm AGL}(1, \mathbb{F}_{2^4})$. We can now extend
this subgroup by an element of order $4$ acting as a field
automorphism of $\mathbb{F}_{2^4}$, giving a unique isomorphism
class of $2$-Frobenius groups of order $2^4\cdot 5\cdot 4$.
Another way of looking at it is that the normalizer of a subgroup
of order $5$ in ${\rm GL}(4, 2)$ is the semilinear group, which
is metacyclic with structure $\mathbb{Z}_{15}: \mathbb{Z}_{4}$,
and this has the Frobenius group $\mathbb{Z}_{5}: \mathbb{Z}_4$
as a subgroup. Reasoning exactly as before, we can show that
there is a unique $2$-Frobenius group of order $3^4\cdot 5\cdot
4$, and it has elementary abelian normal subgroup of order $3^4$.
Now, taking the subdirect product of these gives a unique
isomorphism class of $2$-Frobenius groups of order $2^6\cdot
3^4\cdot 5$. This completes the proof. $\Box$

\noindent{\em Proof of Theorem $\ref{u5(2)}$}. Let $G$ be a finite
group satisfying $${\rm OC}(G)={\rm OC}(U_5(2))=\{2^{10}\cdot
3^5\cdot 5, 11\}.$$  Clearly $|G|=|U_5(2)|=2^{10}\cdot 3^5\cdot
5\cdot 11$ and $s(G)=2$, in fact we have $\pi_1(G)=\{2, 3, 5\}$
and $\pi_2(G)=\{11\}$. Then, by Theorem A in \cite{wili}, one of
the following statements holds:

\begin{itemize}
\item[$(1)$] $G$ is a Frobenius group,

\item[$(2)$]  $G$ is a $2$-Frobenius group, or
\item[$(3)$] $G$ has a normal series $1\unlhd H\lhd K\unlhd G$ such that $H$ is a nilpotent $\pi_1$-group,
$K/H$ is a non-abelian simple group, $G/K$ is a $\pi_1$-group,
and any odd order component of $G$ is equal to one of the odd
order components of $K/H$.
\end{itemize}

If $G$ is a Frobenius group with kernel $K$ and complement $C$,
then ${\rm OC}(G)=\{|K|, |C|\}$, and since $|C|<|K|$, the only
possibility is $|K|=2^{10}\cdot 3^5\cdot 5$ and $|C|=11$.
However, this is a contradiction since $|C|\nmid |K|-1$.

If $G$ is a $2$-Frobenius group of order $2^{10}\cdot 3^5\cdot
5\cdot 11$, then, by the definition, $G=ABC$, where $A$ and $AB$
are normal subgroups of $G$ and $AB$ and $BC$ are Frobenius groups
with kernels $A$ and $B$, respectively. Reasoning as in the proof
of Theorem \ref{u4(2)}, we observe that there are unique
$2$-Frobenius groups $A_1BC$ and $A_2BC$ of orders $2^{10}\cdot
11\cdot 5$ and $3^5\cdot 11\cdot 5$, respectively. Note that,
$A_1$ and $A_2$ are elementary abelian normal subgroups of orders
$2^{10}$ and $3^5$, respectively. Therefore, $G$ is a subdirect
product $(A_1\times A_2)BC=(2^{10}\times 3^5): 11: 5$ of $A_1BC$
and $A_2BC$. So there is a unique $2$-Frobenius group $G=ABC$ of
order $|U_5(2)|=2^{10}\cdot 3^5\cdot 5\cdot 11$.

Finally, we suppose that $G$ satisfies condition (3). Then, by
Table 2, $K/H$ is isomorphic to one of the simple groups
$L_2(11)$, $M_{11}$, $M_{12}$, or $U_5(2)$. We see that, in
general, $K/H\leqslant G/H\leqslant {\rm Aut}(K/H)$. Let $K/H\cong
L_2(11)$. Since $|{\rm Aut}(K/H)|=2^3\cdot 3\cdot 5\cdot 11$ is
not divisible by $3^2$, it follows that $3^4$ divides $|H|$. Let
$P$ be a Sylow $3$-subgroup of $H$ and let $Q$ be a Sylow
$11$-subgroup of $G$ . Then, $P$ is a normal subgroup of $G$,
because $H$ is nilpotent. It now follows that $PQ$ is a subgroup
of $G$ of order $3^4\cdot 11$. Since all groups of order
$3^4\cdot 11$ are nilpotent, we conclude that $3$ is adjacent to
$11$ in ${\rm GK}(G)$, which is a contradiction.

Reasoning exactly as above, we conclude that $K/H\ncong M_{11},
M_{12}$. Therefore, we deduce that $K/H\cong U_5(2)$, and since
$|G|=|U_5(2)|$ it follows that $|H|=1$ and $G=K\cong U_5(2)$.
This completes the proof. $\Box$
\section{Some Open Problems}
We conclude this article with some open problems. Actually, in
this section, we restrict our attention to the relationship
between degree patterns and prime graphs. A natural question is:

\begin{que}\label{question1}
Let $G$ and $M$ be two finite groups with $|G|=|M|$. Clearly
${\rm GK}(G)={\rm GK}(M)$ implies $D(G)=D(M)$. Does the converse
hold?
\end{que}

Assuming the converse is {\em true}, under these hypotheses we
conclude that ${\rm OC}(G)={\rm OC}(M)$, and so $h_{\rm
OD}(M)\leqslant h_{\rm OC}(M)$. In particular, if $M$ is
OC-characterizable, then $M$ is also OD-characterizable. In
\cite[Lemma 2.15]{kogani}, it was shown that if $G$ is a finite
group with $\pi(G)=\pi(M)$ and $D(G)=D(M)$, where $M$ is an
arbitrary alternating or symmetric group, then the prime graphs
of $G$ and $M$ coincide. Therefore, we have the following
consequence.

\begin{corollary}\label{cor-1}
The symmetric and alternating groups which are OC-characterizable
are also OD-characterizable.
\end{corollary}

On the other hand, in view of the Main Theorem in \cite{alavi},
the symmetric groups $S_p$ and $S_{p+1}$, and the alternating
groups $A_p$, $A_{p+1}$ and $A_{p+2}$, where $p\geqslant 3$ is a
prime number, are OC-characterizable. Therefore, by Corollary
\ref{cor-1}, they are also OD-characterizable (see also
\cite[Theorem 1.5]{M.Z(2008)}). We notice that other alternating
and symmetric groups are {\em not} OC-characterizable. It is for
this reason that whose prime graphs are connected and so the
nilpotent groups of the same order have the same order
components. But the situation of OD-characterizability of
alternating and symmetric groups look a little differently. As
pointed out in the Introduction, there are infinitely many
alternating groups $\mathbb{A}_n$ (resp. symmetric groups
$\mathbb{S}_n$) which satisfy $h_{\rm OD}(\mathbb{A}_n)\geqslant
3$ (resp. $h_{\rm OD}(\mathbb{S}_n)\geqslant 3$), in particular,
neither $h_{\rm OD}(\mathbb{A}_n)$ nor $h_{\rm OD}(\mathbb{S}_n)$
is bounded above (see \cite{moghadam}).

We now focus our attention on the sporadic simple groups. By
Table 1 in \cite{degree}, it is easy to see that if $G$ is a
finite group with $\pi(G)=\pi(M)$ and $D(G)=D(M)$, where $M$ is a
sporadic simple group, then the prime graphs of $G$ and $M$
coincide. Moreover, it is proved in \cite{chen} that all sporadic
simple groups are OC-characterizable, hence we conclude that they
are also OD-characterizable (see \cite[Proposition 3.1.]{degree}).

Finally, we consider the OD-characterizability of simple groups of
Lie type. Studies show that between simple groups of Lie type
there are many simple orthogonal and symplectic groups which are
$2$-fold OD-characterizable (see \cite{2-fold}). Moreover, by
Theorem \ref{u4(2)}, we have $h_{\rm OD }(U_4(2))=2$. So far we
have not found a simple group of Lie type $S$ satisfying $h_{\rm
OD}(S)>2$. So it seems natural to ask the following question.
\begin{que}\label{question2}
Does there exist a finite simple group $S$ of Lie type such that
$h_{\rm OD}(S)\geqslant 3$?
\end{que}

\end{document}